\documentclass[12pt]{amsart}
\usepackage{latexsym, amsbsy, amsmath, amsfonts, amssymb, amsthm, amscd, amscd}
\usepackage[all]{xy}
\newcommand{\Z}{\mathbb Z}

\newtheorem{Theorem}{Theorem}[section]

\newtheorem{Corollary}{Corollary}[section]

\allowdisplaybreaks[2]
\numberwithin{equation}{section}
\numberwithin{figure}{section}
\title{On the Lie algebras associated with pure mapping class groups }
\author[Karoui]{R.~Karoui}
\author[Vershinin]{V.~V.~Vershinin}%$^{*}$}
\address{D\'epartement des Sciences Math\'ematiques,
                                     Universit\'e Montpellier II,
Place Eug\`ene Bataillon,
34095 Montpellier cedex 5, France}
\email{Rym.Karoui@math.univ-montp2.fr}
\address{D\'epartement des Sciences Math\'ematiques,
                                     Universit\'e Montpellier II,
Place Eug\`ene Bataillon,
34095 Montpellier cedex 5, France}
\email{ vershini@math.univ-montp2.fr}
\address{Sobolev Institute of Mathematics, Novosibirsk 630090,
Russia } 
\email{ versh@math.nsc.ru}


\begin{document}
\begin{abstract}
Pure braid groups and pure mapping class groups of a punctured sphere
have many features in common.
In the paper the graded Lie algebra of the descending central
series of the pure mapping class of a sphere is studied.
A simple presentation of this Lie algebra is obtained.
\end{abstract}
\maketitle
\tableofcontents

\section{Introduction}

Mapping class group is an  important object in Topology,
Complex Analysis, Algebraic Geometry and other domains. 
It is a lucky case when the method of Algebraic Topology works
perfectly well, the application of the functor of fundamental group
completely solves the topological problem: group of isotopy classes of
homeomorphisms
is described in terms of automorphisms of the fundamental group
of the corresponding surface, as states  the Dehn-Nilsen-Baer theorem, 
see \cite{Iv}, for example. 

Let  $S_{g,b,n}$ be an oriented surface of the genus $g$ with $b$ boundary 
components and with 
a set $Q_n$ of $n$ fixed points. Consider the group 
$\operatorname{Homeo}(S_{g,b,n})$ 
of orientation preserving
self-homeomorphisms of $S_{g,b,n}$ which fix pointwise the boundary (if
it exists) and map the set $Q_n$ into itself. 
Orientation reversing homeomorphisms also possible to consider, see 
\cite{DF}, for example,
but we restrict ourselves to orientation
preserving case. Let 
$\operatorname{Homeo}^0(S_{g,b,n})$ be the normal subgroup of 
self-homeomorphisms of $S_{g,b,n}$ which are isotopic to identity. 
Then the {\it mapping class group} $\mathcal {M}_{g,b,n}$ is defined as a 
factor group

\begin{equation*}
\mathcal {M}_{g,b,n} = \operatorname{Homeo}(S_{g,b,n})/
\operatorname{Homeo}^0(S_{g,b,n}).
\end{equation*}

These groups are connected closely with braid groups.
In \cite{Mag1} W.~Magnus interpreted the braid group as the mapping class
group of a punctured disc with the fixed boundary:
\begin{equation*}
Br_n \cong \mathcal {M}_{0,1,n}. 
\end{equation*}
The same way as the braid groups  the group
$\mathcal {M}_{g,b,n}$ has a natural epimorphism to the symmetric group  
$\Sigma_{n}$ with the kernel called the  {\it pure mapping class group} $\mathcal {PM}_{g,b,n}$,
so there exists an exact sequence:
\begin{equation*}
1 \to \mathcal {PM}_{g,b,n} \to \mathcal {M}_{g,b,n} \to \Sigma_{n}\to 1.
\end{equation*}
Geometrically  pure mapping class group $\mathcal {PM}_{g,b,n}$ 
is descried as consisting of isotopy classes of homeomorphisms that preserve the  punctures pointwise.

We call an element $m$ of the pure mapping class group 
$ \mathcal{PM}_{g,b,n}$ $i$-{\it Makanin} or $i$-{\it Brunnian} if 
a homeomorphism $h$ lying in the class $m$
\begin{equation*}
h: S_{g,b,n} \to S_{g,b,n},
\end{equation*}
becomes  isotopical to the identity map if we fill the deleted point $i$
(or if it becomes non-fixed). Filling the point $i$ generates the homomorphism
$$pm_i: \mathcal {PM}_{g,b,n} \to \mathcal {PM}_{g,b,n-1}.$$  
 We  denote the subgroup
of $i$\,-Makanin elements of the mapping class group by $A_i$,
it is so the kernel of $pm_i$. The subgroups $A_i$, $i=1, \dots, n$,
 are conjugate in $\mathcal {M}_{g,b,n} $. 
The intersection of all subgroups of $i$\,-Makanin elements is the of 
{\it Makanin} or {\it Brunnian} subgroup of the mapping class group
\begin{equation*} 
Mak_{g,b,n} =\cap_ {i=1}^{n} A_i.
\label{eq:mak}
\end{equation*}

In the paper we consider the pure mapping class group of a sphere with
no boundary components $\mathcal {PM}_{0,0,n}$, which we denote for
simplicity by $\mathcal {PM}_{0,n}$. We study the  
 natural Lie algebra
obtained from the descending central series for $\mathcal {PM}_{0,n}$. 
One motivation
for the work here is that the group $\mathcal {PM}_{0,n}$ 
is natural as well as accessible case continuing the same study for
the pure braid group done in the works of T.~Kohno \cite{K}, 
T.~Kohno and T.~Oda \cite{KO}, Y.~Ihara \cite{Ih}, R.~Bezrukavnikov \cite{Bez}
and for McCool group it was done in the works \cite{CPVW}
and \cite{BerPa}.

%The techniques here for addressing these Lie
%algebras are due to T.~Kohno \cite{K} and M.~Falk and R.~Randell
%\cite{FR}.

\section{Lie algebra $gr^*(\mathcal {PM}_{0,n})$ %\label{sec:ib_imc}
}

The group $\mathcal {PM}_{0,n}$ is closely related to the pure braid group
on a sphere $P_n(S^2)$ as well as its non-pure analogue $\mathcal {M}_{0,n}$
is connected with the (total) braid group of a sphere $Br_n(S^2)$. 

We start with presentations.
Usually  the braid group $Br_n$ 
is given by the following Artin presentation \cite{Art1}.
It has the generators $\sigma_i$, 
$i=1, ..., n-1$, and two types of relations: 
\begin{equation}
 \begin{cases} \sigma_i \sigma_j &=\sigma_j \, \sigma_i, \ \
\text{if} \ \ |i-j| >1,
\\ \sigma_i \sigma_{i+1} \sigma_i &= \sigma_{i+1} \sigma_i \sigma_{i+1} \, .
\end{cases} \label{eq:brelations}
\end{equation}

The generators $a_{i,j}$, $1\leq i<j\leq n $ for the pure
braid group $P_n$ (of a disc) can be defined (as elements of the the braid 
group $Br_n$)  by the formula:
$$a_{i,j}=\sigma_{j-1}...\sigma_{i+1}\sigma_{i}^2\sigma_{i+1}^{-1}...
\sigma_{j-1}^{-1}.$$ 
Then the defying relations, which are called the \emph{Burau
relations} \cite{Bu1}, \cite{Mar2} are as follows: 
%\begin{multline}
\begin{equation}
\begin{cases}
a_{i,j}a_{k,l}=a_{k,l}a_{i,j}
\ \text {for} \ i<j<k<l \ \text {and} \ i<k<l<j, \\
a_{i,j}a_{i,k}a_{j,k}=a_{i,k}a_{j,k}a_{i,j} \ \text {for} \
i<j<k, \\
a_{i,k}a_{j,k}a_{i,j}=a_{j,k}a_{i,j}a_{i,k} \ \text
{for} \ i<j<k, \\
a_{i,k}a_{j,k}a_{j,l}a_{j,k}^{-1}=a_{j,k}a_{j,l}a_{j,k}^{-1}a_{i,k}
\ \text {for} \ i<j<k<l.\\
\end{cases}
\label{eq:burau}
\end{equation}
%{multline}

\smallskip

\noindent
It was proved by O.~Zariski \cite{Za1} and then
rediscovered by E.~Fadell and J.~Van Buskirk \cite{FaV} that a presentation 
for the braid group
of a sphere can be given with the   generators 
$\sigma_i$, $i=1, ..., n-1$, the same as for the classical braid group, satisfying
the braid relations (\ref{eq:brelations})
and the following sphere relation: 
\begin{equation}
\sigma_1 \sigma_2 \dots \sigma_{n-2}\sigma_{n-1}^2\sigma_{n-2} \dots
\sigma_2\sigma_1 =1.
\label{eq:spherelation}
\end{equation}
Having the same generators, but if we add to the braid relations
(\ref{eq:brelations}) and the sphere relation 
(\ref{eq:spherelation})
one more relation 
((\ref{eq:sphe_mc}) below)
we get the presentation for the 
mapping class group of a punctured sphere $\mathcal {M}_{0,n}$
 obtained by W.~Magnus \cite{Mag1}, see also \cite{Mag3} and 
 \cite{MKS}. 
\begin{equation}
(\sigma_1 \sigma_2 \dots \sigma_{n-2}\sigma_{n-1})^n =1.
\label{eq:sphe_mc}
\end{equation}
 
Let $\Delta$ be the Garside's   fundamental word   in the braid 
group $Br_{n}$ \cite{Gar}. It can be in particular expressed by the formula:
$$\Delta = \sigma_1 \dots \sigma_{n-1} \sigma_1 \dots \sigma_{n-2} \dots  
\sigma_1 \sigma_2 \sigma_1.$$
If we use Garside's notation $\Pi_t\equiv \sigma_1\dots \sigma_t$, then
$\Delta \equiv \Pi_{n-1} \dots \Pi_1$.
If the generators $\sigma_1,$ $\sigma_2$, $\dots$, $\sigma_{n-2},$
$\sigma_{n-1}$,
are subject to the braid relations (\ref{eq:brelations}), then the 
condition (\ref{eq:sphe_mc}) is equivalent to the following relation
\begin{equation*}
\Delta^2 =1.
\label{eq:sphe_mcD}
\end{equation*}
For the pure braid group on a sphere let us introduce the elements
$a_{i,j}$ for all $i, j$ by the formulas:  
%\begin{multline}
\begin{equation}
\begin{cases}
a_{j,i}= a_{i,j} \ \ \text{for} \ i<j\leq n,\\
a_{i,i}= 1. 
\label{eq:aji}
%\end{equation*}
\end{cases}
\end{equation}
%{multline}
The pure braid group for the sphere has the generators $a_{i,j}$
which satisfy Burau relations (\ref{eq:burau}), relations (\ref{eq:aji}),
 and the following relations \cite{GVB}:
\begin{equation*}
a_{i,i+1}a_{i,i+2} \dots a_{i,i+n-1} = 1 \ \ \text{for all} \ i\leq n,\\
\end{equation*}
with the convention that $k+n =k$.
Note that $\Delta^2$ is a pure braid and it can be expressed by the following
formula
%\begin{equation}
%$$
\begin{multline*}
\Delta^2 = (a_{1,2}a_{1,3} \dots a_{1,n})(a_{2,3}a_{2,4} \dots a_{2,n})
\dots (a_{n-1,n}) = \\
(a_{1,2})(a_{1,3} a_{2,3})( a_{1,4}a_{2,4}a_{3,4}) \dots (a_{1,n} \dots
\dots a_{n-1,n}).
\end{multline*}
%\end{equation}
%$$
The fact that this element of the braid group generates its center goes back to
Chow \cite{Ch}. 

Let us denote by $P_n(S^2_3)$ the pure braid group on $n$ strings of a sphere 
with three points deleted or equivalently the subgroup of the pure braid group 
of a disc on $n+2$ strings where (say, the last) two strings are fixed. 

The following statement follows from the normal forms of the groups
$P_n(S^2)$ and $\mathcal {PM}_{0,n}$ \cite{GVB} and on the geometrical level 
it 
was expressed in \cite{GG2}. Note that the groups $P_2(S^2)$ and $\mathcal
 {PM}_{0,3}$ are trivial.
\begin{Theorem} (i) The pure braid group of a sphere $P_n(S^2)$ for $n\geq 3$
is isomorphic to the direct product
of the cyclic group of order 2 (generated by $\Delta^2$) and $\mathcal {PM}_{0,n}$.
\par
(ii) The pure braid group $P_n$ for $n\geq 2$ is isomorphic to the direct product
of the infinite cyclic group (generated by $\Delta^2$) and $\mathcal {PM}_{0,n+1}$. 
\par 
(iii) The groups $\mathcal {PM}_{0,n}$ and $P_{n-3}(S^2_3)$ are isomorphic for 
$n\geq 4$.
\label{theo:gil-vb} 
\end{Theorem} 

The isomorphism of the part (i) of Theorem~\ref{theo:gil-vb} is compatible 
with the homomorphisms
$p_i:P_n(S^2)\to P_{n-1}(S^2)$,  
$pm_i : \mathcal{PM}_{0,n} \to \mathcal{PM}_{0,n-1}$  
consisting of deleting one string or forgetting one point, so the group of 
Makanin braids of a
 sphere coincide with the subgroup of Makanin mapping class of a sphere.

For a group $G$ the descending central series 
\begin{equation*}
G =\Gamma_1  > \Gamma_2 > \dots  > \Gamma_i > \Gamma_{i+1} > \dots .
\end{equation*}
\noindent
is defined by the formulas
\begin{equation*}
\Gamma_1 = G, \ \ \Gamma_{i+1} =[\Gamma_{i}, G].
\end{equation*}

The descending central series of a discrete group $G$ gives rise to the 
associated graded Lie algebra (over $\Z$) $gr^*(G)$ \cite{Serr}. 
\begin{equation*}
gr^i(G)= \Gamma_i/\Gamma_{i+1}.
\end{equation*}

%The descending central series of a discrete group $G$ gives rise to the 
%associated graded Lie algebra (over $\Z$) $gr^*(G)$ \cite{Serr}.  

A presentation of the Lie algebra $gr^*(P_n)$ for the pure braid group
can be described as follows \cite{K}. It is the quotient of the
free Lie algebra $L[A_{i,j}| \, 1 \leq i < j \leq n]$ generated by
elements $A_{i,j}$ with $1 \leq i < j \leq n$ modulo the
``infinitesimal braid relations" or ``horizontal $4T$ relations"
given by the following three relations:

\begin{equation*}
\begin{cases}
 [A_{i,j}, A_{s,t}] = 0, \  \text{if} \ \{i,j\} \cap \{s,t\} = \phi, \\
 [A_{i,j}, A_{i,k} + A_{j,k}] = 0, \  \text{if} \ i<j<k , \\
  [A_{i,k}, A_{i,j} + A_{j,k}] = 0, \ \text{if} \ i<j<k. \\
  \end{cases}
  \label{eq:kohno}
\end{equation*}
Y.~Ihara in \cite{Ih} gave  a presentation of the Lie algebra $gr^*(P_n(S^2))$ 
of the pure braid group of a sphere. It is convenient to 
have conventions like (\ref{eq:aji}). So, 
it is the quotient of the
free Lie algebra $L[B_{i,j}| \, 1 \leq i ,j \leq n]$ generated by
elements $B_{i,j}$ with $1 \leq i , j \leq n$ modulo the
 following  relations:
\begin{equation*}
\begin{cases}
B_{i,j} =  B_{j,i} \  \text{for} \ 1\leq i,j \leq n , \\
B_{i,i} =  0 \ \ \text{for} \ 1\leq i \leq n , \\
 [B_{i,j}, B_{s,t}] = 0, \  \text{if} \ \{i,j\} \cap \{s,t\} = \phi, \\
 \sum_{j=1}^n B_{i,j} = 0, \  \text{for} \ 1\leq i \leq n. \\
    \end{cases}
    \label{eq:ihara}
\end{equation*}

\smallskip

\noindent
%It is interresting to note that the Kohno relations (type three and four in 
%(\ref{eq:kohno})) are the consequences of third and the fourth type relations
%in (\ref{eq:ihara}).
It is a factor algebra of the algebra $gr^*(P_n)$: the last two relations in 
(\ref{eq:kohno}) are the consequences of the third and the forth type 
relations in (\ref{eq:ihara}).

\begin{Theorem} (i) The graded Lie algebra $gr^*(\mathcal {PM}_{0,n})$ 
is the quotient of the
free Lie algebra $L[B_{i,j}| \, 1 \leq i ,j \leq n]$ modulo the
 following  relations:
\begin{equation}
\begin{cases}
B_{i,j} =  B_{j,i} \  \text{for} \ 1\leq i,j \leq n , \\
B_{i,i} =  0 \ \ \text{for} \ 1\leq i \leq n , \\
 [B_{i,j}, B_{s,t}] = 0, \  \text{if} \ \{i,j\} \cap \{s,t\} = \phi, \\
 \sum_{j=1}^n B_{i,j} = 0, \  \text{for} \ 1\leq i \leq n, \\
 \medskip
 \sum_{i=1}^{n-1}\sum_{j=i+1}^n B_{i,j} = 0. \\ 
    \end{cases}
    \label{eq:m0n}
\end{equation} 
(ii) The graded Lie algebra $gr^*(\mathcal {PM}_{0,n})$ is the
quotient of the free Lie algebra
$L[A_{i,j}| \, 1 \leq i < j \leq n-1]$ generated by
elements $A_{i,j}$ with $1 \leq i < j \leq n-1$ modulo the
 following  relations:
\begin{equation*}
\begin{cases}
 [A_{i,j}, A_{s,t}] = 0, \  \text{if} \ \{i,j\} \cap \{s,t\} = \phi, \\
 \smallskip
 \sum_{i=1}^{n-2}\sum_{j=i+1}^{n-1} A_{i,j} = 0. \\ 
   \end{cases}
  \label{eq:m0n2}
\end{equation*}
\label{theo:lie_m0n} 
\end{Theorem} 
\begin{proof} 
Part (i) of Theorem~\ref{eq:m0n} follows from the Ihara presentation and the 
part (i) of Theorem~\ref{theo:gil-vb}. The pure mapping class group  
$\mathcal {PM}_{0,n}$ 
is a direct summand in the pure braid group of a sphere $P_n(S^2)$,
so there is no  problem in obtaining exact sequence after application
of the functor of associated graded Lie algebras $gr^*$. 

To obtain part (ii) let us write in detail the system of linear equations 
which constitute the forth part of the Ihara relations (\ref{eq:ihara}):
\begin{equation}
\begin{cases}
 B_{1,2} + B_{1,3} + \dots B_{1,n} = 0, \\  
  B_{1,2} + B_{2,3} + \dots B_{2,n} = 0, \\ 
  \dots \\
   B_{1,n} + B_{2,n} + \dots B_{n-1,n} = 0. \\ 
\end{cases}
\label{eq:4lie_sph}
\end{equation}
\noindent
The $n-1$ equations (except the last one) give the possibility to exclude
the letters $B_{i,n} \  \text{for} \ 1\leq i \leq n-1$ from the presentation.
Note that the third type relations in (\ref{eq:ihara}) with $j=n$ 
are the consequences of the same type relation with $j\leq n-1$ and  
 relations (\ref{eq:4lie_sph}) except the last one. Take then the 
 linear combination of these  equations where the first $n-1$
equations are taken with the coefficient $+1$ and the last one  with the
coefficient $-1$. We get the equation 
\begin{equation*}
2 (\sum_{i=1}^{n-2}\sum_{j=i+1}^{n-1} A_{i,j}) = 0. \\ 
  %  \label{eq:kohno}
\end{equation*}
The second type relation in (\ref{eq:m0n}) is a consequence of the 
last equation in (\ref{eq:ihara}) and first relation in (\ref{eq:4lie_sph}).
\end{proof}
\begin{Corollary}
A presentation of the Lie algebra $gr^*(P_n(S^2))$ can be given with generators
$A_{i,j}$ with $1 \leq i< j  \leq n-1$,  modulo 
%the relations (\ref{eq:m0n2}) and 
the following  relations:
\begin{equation*}
\begin{cases}
 [A_{i,j}, A_{s,t}] = 0, \  \text{if} \ \ \{i,j\} \cap \{s,t\} = \phi, \\
 \smallskip
2  (\sum_{i=1}^{n-2}\sum_{j=i+1}^{n-1} A_{i,j}) = 0. 
 \end{cases}
  %  \label{eq:kohno}
\end{equation*}
So, the element $\sum_{i=1}^{n-2}\sum_{j=i+1}^{n-1} A_{i,j}$ of order 2 generates
the central subalgebra in  $gr^*(P_n(S^2))$. 
\end{Corollary}
\hfill $\square$

\section{Example\label{sec:ex}}

The pure braid group of a sphere $P_4(S^2)$  is isomorphic to the direct product
of the cyclic group of order 2 (generated by $\Delta^2$) and the pure braid group on one string of a sphere with three points deleted, that is the fundamental group 
of disc with two points deleted, that is a free group on two generators $F_2$.
Its associated graded Lie algebra is a direct sum of central $\Z/2$ and the
free Lie algebra on two generators. The pure mapping class group 
$\mathcal{PM}_{0,4}$ is isomorphic to a free group on two generators. According 
to Theorem~\ref{theo:lie_m0n} its associated graded Lie algebra is the
quotient of the free Lie algebra
$L[A_{1,2}, A_{1,3}, A_{2,3}]$ modulo the 
 following  relation:
\begin{equation*}
A_{1,2} + A_{1,3} + A_{2,3}  = 0, 
\end{equation*}
so, is a free Lie algebra on two generators.


\begin{thebibliography}{References}
\bibitem{Art1}
\emph{ E. Artin},
 Theorie der Z\"opfe. 
 Abh. Math. Semin. Univ. Hamburg, 1925,
v. 4, 47--72.
%\bibitem{BaM}
%\emph{V.~Bardakov; R.~Mikhailov}, On certain questions of the free group %automorphisms theory.  Comm. Algebra  36  (2008),  no. 4, 1489--1499.
\bibitem{BerPa}
\emph{B.~Berceanu, S.~Papadima},
 Universal representations of braid and braid-permutation groups, arXiv:0708.0634 
\bibitem{BCWW}
\emph{J.~A.~Berrick, F.~R.~Cohen, Y.~L.~Wong and J.~Wu}, Configurations, 
braids, and homotopy groups.  J. Amer. Math. Soc.  19  (2006),  no. 2, 265--326.
\bibitem{Bez}
\emph{R.~Bezrukavnikov}, Koszul DG-algebras arising from configuration spaces.  
Geom. Funct. Anal.  4  (1994),  no. 2, 119--135.
 \bibitem{Bu1}
\emph{W.~Burau}, 
\"Uber Zopfinvarianten. (German)
Abh. Math. Semin. Hamb. Univ. 9, 117-124 (1932). 
\bibitem{Ch}
\emph{W.-L. Chow}, On the algebraical braid group,
 Ann. Math. 1948,  49{\rm , No 3}, 654--658.
\bibitem{CPVW} \emph{F.~R.~Cohen;  J.~Pakianathan; V.~V.~Vershinin; J.~Wu}, 
Basis-conjugating automorphisms of a free group and associated Lie algebras. 
Iwase, Norio (ed.) et al., Proc. of the conference on groups, homotopy and configuration spaces, Univ. of Tokyo, July, 2005. Coventry: Geometry and Topology Monographs 13, 147-168 (2008).
\bibitem{DF}
\emph{W.~Dicks and E.~Formanek},
Algebraic mapping-class groups of orientable surfaces with boundaries.
pp. 57-115, in:   Infinite groups: geometric, combinatorial and dynamical aspects.
Progress in Math. 248,   Birkh\"auser Verlag, Basel, 2005.
\bibitem{FaV}
\emph{E.~Fadell and J.~Van Buskirk}, The braid groups of $E\sp{2}$ and $S\sp{2}$. 
Duke Math. J. 29 1962,    243--257.
\bibitem{FR} \emph{M.~Falk, and R.~Randell}, The lower central series of a
fiber--type arrangement, Invent. Math., {\bf{82}} (1985), 77--88.
\bibitem{FLM}
\emph{B.~Farb, C.~J.~Leininger, D.~Margalit} The lower central series and pseudo-Anosov dilatations. arXiv:math/0603675.  
\bibitem{Gar}
\emph{F.~A.~Garside}, 
The braid group and other groups,
Quart. J. Math. Oxford Ser. 1969, 20,
235--254.
\bibitem{GVB}
\emph{R.~Gillette; J.~Van Buskirk}, The word problem and consequences for 
the braid groups and mapping class groups of the $2$-sphere.  
Trans. Amer. Math. Soc.  131.  1968. 277--296. 
\bibitem{GG2}
\emph{D.~L.~Gonçalves; J.~Guaschi}, The roots of the full twist for surface 
braid groups. Math. Proc. Cambridge Philos. Soc. 137 (2004), no. 2, 307--320. 
\bibitem{Iv}
\emph{N. V. Ivanov}, Mapping class groups, Handbook of geometric topology, 523--633, North-Holland, Amsterdam, 2002. 
\bibitem{Ih}
 \emph{Y.~Ihara},   Galois group and some arithmetic functions,
Proceedings of the International Congress of mathematicians, Kyoto,
1990, Springer (1991),  99-120. 
\bibitem{K} \emph{T.~Kohno}, S\'erie de Poincar\'e-Koszul associ\'ee aux
groupes de tresses pure, Invent. Math., {\bf{82}} (1985), 57-75.
\bibitem{KO} 
\emph{T.~Kohno; T.~Oda}, The lower central series of the pure braid group of an algebraic curve.  Galois representations and arithmetic algebraic geometry (Kyoto, 1985/Tokyo, 1986),  201--219, Adv. Stud. Pure Math., 12, North-Holland, Amsterdam, 1987. 
\bibitem{Mag1}
\emph{W.~Magnus},
\"Uber Automorphismen von Fundamentalgruppen berandeter Fl\"achen. (German)
 Math. Ann. 109, 617-646 (1934). 
%\bibitem{Mag2}
%\emph{W.~Magnus},
% Braids and Riemann surfaces. Comm. Pure Appl. Math. 25 (1972) 151-161.
 \bibitem{Mag3}
 \emph{W.~Magnus},
 Braid groups: A survey.  Proceedings of the Second International Conference on the Theory of Groups (Australian Nat. Univ., Canberra, 1973),  pp. 463--487. Lecture Notes in Math., Vol. 372, Springer, Berlin, 1974. 
\bibitem{MKS}
\emph{W.~Magnus; A.~Karrass; D.~Solitar},
Combinatorial group theory. Presentations of groups in terms of generators 
and relations. 2nd rev. ed. 
Dover Books on Advanced Mathematics. New York: Dover Publications, 
Inc. XII, 444 p. (1976). 
\bibitem{Mar2}
\emph{A. A. Markoff},  Foundations of the Algebraic Theory of
Tresses, Trudy Mat. Inst. Steklova, No~16, 1945 (Russian, English
summary).
\bibitem{Sc}
\emph{G.~P.~Scott}, Braid groups and the group of homeomorphisms of a 
surface. Proc. Cambridge Philos. Soc. 68, 1970,  605--617.
\bibitem{Serr}
\emph{J.-P. Serre},
 Lie algebras and Lie groups. 1964 lectures given at Harvard University. Corrected fifth printing of the second (1992) edition. Lecture Notes in Mathematics, 1500. Springer-Verlag, Berlin, 2006. viii+168 pp.
\bibitem{Ve10}
\emph{V.~V.~Vershinin}, Braid groups, their Properties and Generalizations.
Handbook of Algebra, vol.~4, Elzseier, Amsterdam a.o. 2006, p.~427-465.  
\bibitem{Za1}
\emph{O.~Zariski}, On the Poincare
group of rational plane curves,
Am. J. Math.
 1936, 58,  607-619. 
\end{thebibliography}
\end{document}